\documentclass[amsfonts,11pt]{amsart}
\RequirePackage{amsmath,amssymb,amsthm, amscd, comment,mathtools}
\usepackage[margin=1.2in]{geometry}
\usepackage{amsmath}
\usepackage{amsthm}
\usepackage{amssymb}
\usepackage{amscd}
\usepackage{mathrsfs}
\usepackage{graphicx}
\usepackage{rotating, graphicx}
\usepackage[all,cmtip]{xy}
\usepackage{enumerate}
\usepackage{tikz-cd}
\usepackage[thinlines]{easytable}
\usepackage{multirow}
\usepackage{mathtools}
\usepackage{mathrsfs}
\usepackage{changepage} % for \adjustwidth
\usepackage{paralist}%\usepackage{enumerate}
\usepackage{MnSymbol}
\usepackage{kantlipsum}
\usepackage{enumitem}
\usepackage{times}
\usepackage{subfig}
\usepackage{color}
\usepackage[pagebackref,breaklinks,colorlinks,linkcolor=red,anchorcolor=red,citecolor=blue]{hyperref}

\calclayout

\newtheorem{proposition}[subsubsection]{Proposition}
\newtheorem{corollary}[subsubsection]{Corollary}
\newtheorem{theorem}[subsubsection]{Theorem}
\newtheorem{lemma}[subsubsection]{Lemma}
\theoremstyle{definition}
\newtheorem{definition}[subsubsection]{Definition}
\theoremstyle{remark}

\newtheorem{remark}[subsubsection]{Remark}

\numberwithin{equation}{subsubsection}

\DeclareMathAlphabet{\mathbbold}{U}{bbold}{m}{n}

\title{Some results on the motivic nearby cycle}

\author{Fangzhou Jin}
\address{School of Mathematical Sciences\\
Tongji University\\
Siping Road 1239\\
200092 Shanghai\\
P.R.China}
\email{\href{mailto:fangzhoujin@tongji.edu.cn}{fangzhoujin@tongji.edu.cn}}
\urladdr{\url{https://fangzhoujin.github.io/}}

\author{Enlin Yang}
\address{School of Mathematical Sciences\\
Peking University\\
No.5 Yiheyuan Road Haidian District,\\
Beijing 100871,\\
P.R.China}
\email{\href{mailto: yangenlin@math.pku.edu.cn}{yangenlin@math.pku.edu.cn}}
\urladdr{\url{https://www.math.pku.edu.cn/teachers/yangenlin/ely.htm}}

\date{\number\day-\number\month-\number\year}

\begin{document}

\maketitle

\begin{abstract}

We extend Ayoub's formalism of motivic nearby cycle functor to the $\infty$-categorical level, and prove some desired cohomological properties by relating the motivic nearby cycle functor to the notion of local acyclicity in motivic homotopy.
\end{abstract}

\tableofcontents

\noindent

\section{Introduction}
\label{sec:nearby}

\subsection{Towards the motivic nearby cycle}

\subsubsection{}
The formalism of nearby cycles for \'etale sheaves is constructed by Grothendieck by analogy with the \emph{Milnor fibration} in complex geometry (\cite[XIII.1.3.2]{SGA7}): if $S$ is the spectrum of a henselian discrete valuation ring with generic point $\eta$ and closed point $\sigma$, and if $X$ is an $S$-scheme, then the \emph{nearby cycle functor} is defined as
\begin{align}
\begin{split}
R\Psi_X:D^+(X_\eta,\Lambda)&\to D^+(X_{\bar{\sigma}},\Lambda)\\
\mathcal{F}&\mapsto \bar{i}_X^*R\bar{j}_{X*}(\mathcal{F}_{\bar{\eta}}).
\end{split}
\end{align}
Here $\bar{\eta}$ is a geometric generic point defining the geometric point $\bar{\sigma}$ over $\sigma$, and the diagram $X_{\bar{\sigma}}\xrightarrow{\bar{i}_X}X_{(\bar{\sigma})}\xleftarrow{\bar{j}_X}X_{\bar{\eta}}$ is obtained from the diagram $\bar{\sigma}\xrightarrow{}S_{(\bar{\sigma})}\xleftarrow{}\bar{\eta}$ by base change, where $S_{(\bar{\sigma})}$ is the strict localization of $S$ at $\bar{\sigma}$. Via such a construction, Grothendieck establishes the \emph{geometric local monodromy theorem} (\cite[I Th. 1.2]{SGA7}), solving and vastly generalizing a conjecture of Milnor on the quasi-unipotence of the monodromy action.

\subsubsection{}
In motivic homotopy theory, given the six functor formalism (\cite{Ayo}, \cite{CD}), a natural analogue for such a functor in the setting above appears to be the following one:
\begin{align}
\begin{split}
\Psi'_X:\mathcal{SH}(X_\eta)&\to \mathcal{SH}(X_{\bar{\sigma}})\\
\mathcal{F}&\mapsto \bar{i}_X^*R\bar{j}_{X*}(\mathcal{F}_{\bar{\eta}}).
\end{split}
\end{align}
As Ayoub points out (\cite[\S1.4]{Ayo1}), one of the major drawbacks of this functor lies in the fact that we work with integral coefficients instead of torsion ones, and as a consequence one has $\Psi'_S(\mathbb{Z})\neq\mathbb{Z}$ in general.  Instead, in his thesis (\cite[\S 3]{Ayo}), Ayoub takes a different path by constructing a cosimplicial motive $\mathcal{A}^\bullet$ which is a motivic analogue of the Rapoport-Zink complex (\cite{RZ}), and defines the \emph{unipotent motivic nearby cycle functor} $\Upsilon$ by taking a modification of the \emph{canonical specialization system} (see \ref{num:cansp}) by $\mathcal{A}^\bullet$. By further taking into account the finite \'etale covers, he constructs the \emph{(tame) motivic nearby cycle functor} $\Psi$ out of $\Upsilon$. We will recall his construction in Section~\ref{sec:ayocons}. Ayoub's main results on $\Psi$ state that in characteristic $0$, it preserves constructible objects, commutes with duality, and satisfies a K\"unneth formula (\cite[Th. 3.5.14, 3.5.17, 3.5.20]{Ayo}). For further applications of this formalism, see for example \cite{Ayo1}, \cite{AIS}, \cite{LPS}, \cite{Azo} and \cite{Ayo3}.

\subsection{Summary of results}

\subsubsection{}
In \cite[Cor. 3.5.18]{Ayo}, Ayoub proves that in characteristic $0$, the tame motivic nearby cycle functor $\Psi$ satisfies a K\"unneth formula, that is, the formation of $\Psi$ commutes with exterior products over the base scheme $\mathbb{A}^1$. This property fails in positive characteristic, already for \'etale cohomology (see \cite[Rem. 2.2]{Ill}). On the other hand, working with \'etale motives, Ayoub proves a K\"unneth formula for the \emph{total motivic nearby cycle functor} instead of the tame one (\cite[Th. 10.19]{Ayo2}).

\subsubsection{}
Motivated by the K\"unneth formula in the \'etale setting (\cite[Th. A.3]{Ill}), we use a different approach here: we prove a K\"unneth formula for the tame motivic nearby cycle $\Psi$ while restricting to certain motives with nice properties, given by the notion of local acyclicity in the motivic world that we introduced in \cite[Def. 2.1.7]{JY} (see Definition~\ref{def:locacy} below). 
In \'etale cohomology, it is well-known that nearby cycles can be used to detect local acyclicity (see Remark~\ref{rm:lavan} below); for motives, we prove a result following this spirit (Proposition~\ref{prop:lakun}): if a motive over the total space is locally acyclic over the base, then its motivic nearby cycle is simply the restriction to the special fiber. We then deduce a K\"unneth formula (Theorem~\ref{th:lakun}) which states that $\Psi$ commutes with with exterior products when restricted to objects that are locally acyclic over the base. In particular, we prove the monoidality of $\Psi$ when applied to globally dualizable objects (Corollary~\ref{cor:psimon}). Note that these results do not require resolution of singularities or inverting $p$.

\subsubsection{}
Here is the organization of the results of the paper. 
In Section~\ref{sec:sps}, we develop a notion of specialization systems in the $\infty$-categorical language.
We enhance Ayoub's construction of the motivic nearby cycle functors to $\infty$-functors in Section~\ref{sec:ayocons}.

\subsubsection{}
The upshot of Section~\ref{sec:kun} is to establish a K\"unneth formula for $\Psi$. We first recall the notion of local acyclicity in motivic $\infty$-categories (Definition~\ref{def:locacy}), and show that this notion is stable under exterior products (Lemma~\ref{lm:laprod}). In Proposition~\ref{prop:lakun}, we prove that the functor $\Psi j^*$ applied to locally acyclic objects agrees with the restriction to the special fiber. From these properties we deduce the following K\"unneth formula:
\begin{theorem}[Theorem~\ref{th:lakun}]
Let $X$ and $Y$ be two $\mathbb{A}^1_k$-schemes and let $K\in \mathcal{SH}(X)$ and $L\in \mathcal{SH}(Y)$ be two objects that are universally strongly locally acyclic (USLA) over $\mathbb{A}^1_k$ in the sense of Definition~\ref{def:locacy}, then there is a canonical isomorphism
\begin{align}
\Psi_X(K_\eta)\boxtimes_{\mathbb{A}^1_k}\Psi_Y(L_\eta)
\to
\Psi_{X\times_{\mathbb{A}^1_k}Y}(K\boxtimes_{\mathbb{A}^1_k}L)_\eta.
\end{align}
\end{theorem}
Note that this result holds without inverting $p$ or assuming resolution of singularities. In particular we deduce the following monoidality of $\Psi$:
\begin{corollary}[Corollary~\ref{cor:psimon}]
The functor $\Psi_{Id}j^*:\mathcal{SH}(\mathbb{A}^1_k)\to\mathcal{SH}(\sigma)$ is symmetric monoidal when restricted to the full subcategory of objects that are USLA relatively to the identity of $\mathbb{A}^1_k$.
\end{corollary}
For example, any dualizable object in $\mathcal{SH}(\mathbb{A}^1_k)$ satisfies this property. The converse is also true under resolution of singularities, see Proposition~\ref{prop:ladual}. In Proposition~\ref{prop:labc} we show that $\Psi$ commutes with locally acyclic base changes, including for example smooth morphisms or any morphism of the form $X\to k$ assuming resolution of singularities.

\subsubsection*{\bf Acknowledgments}
We would like to thank Joseph Ayoub for pointing out an error in a previous version. Both authors are supported by the National Key Research and Development Program of China Grant Nr.2021YFA1001400. F. Jin is supported by the National Natural Science Foundation of China Grant Nr.12101455 and the Fundamental Research Funds for the Central Universities. E. Yang is supported by Peking University's Starting Grant Nr.7101302006 and NSFC Grant Nr.11901008.

\subsection{Notations and conventions}

\subsubsection{}
We denote by $\Delta$ the category of simplexes, or equivalently the category of finite linearly ordered sets with (non-strictly) order-preserving maps. We denote by $CAlg(Pr_{Stab})$ the $\infty$-category of symmetric monoidal presentable stable $\infty$-categories (\cite{HA}).

\subsubsection{}
All schemes are noetherian (in fact most of the time we work with schemes of finite type over a field of positive characteristic). Smooth morphisms are assumed to be separated of finite type.

\subsubsection{}
We use the language of motivic $\infty$-catgegories (\cite{Kha}). An important example is the motivic stable homotopy category $\mathcal{SH}$, which turns out to be the universal such category (\cite[Th. 7.14]{DG}).

Apart from the six functors formalism, the definition of the motivic nearby cycle functor makes use of $(\infty,1)$-categorical colimits, which corresponds to homotopy colimits in the language of model categories.

\subsubsection{}
\label{num:notB}
Let $k$ be a field of characteristic $p$. Let $B$ be a smooth connected $1$-dimensional $k$-scheme and let $\pi$ be a global section of $B$. We denote by $i:\sigma=B/(\pi)\to B$ the closed immersion with $j:\eta=B-\sigma\to B$ the open complement. For any $B$-scheme $X$, we denote by $j_X:X_\eta\to X$ and $i_X:X_\sigma\to X$ the base change of $j$ and $i$.

\subsubsection{}
If $f:X\to Y$ is a morphism and $K$ is an object in $\mathcal{T}(Y)$, we denote $K_{|X}=f^*K\in\mathcal{T}(X)$.

\subsubsection{}
Let $S$ be a scheme and let $X_1$ and $X_2$ be two $S$-schemes. Denote by $p_i:X_1\times_SX_2\to X_i$ the projections, $i=1,2$.
For $K_1\in {\mathbf T}(X_1)$ and $K_2\in {\mathbf T}(X_2)$, we denote
\begin{align}
K_1\boxtimes_S K_2\coloneqq p_1^* K_1\otimes p_2^* K_2 
\end{align}
called the \textbf{exterior product}, which is an object of $\mathcal{T}(X_1\times_S X_2)$.

\subsubsection{}
Recall that a motivic $\infty$-catgegory $\mathcal{T}$ is \textbf{continuous} if for every cofiltered system $(S_\alpha)_{\alpha}$ of schemes with affine transition maps and with limit $S$, the canonical functor
\begin{align}
\operatorname{colim}\mathcal{T}(S_\alpha)\to \mathcal{T}(S)
\end{align}
is an equivalence (\cite[Def. A.2]{DFJK}). Here the colimit is taken in the $\infty$-category of presentable $\infty$-categories and colimit-preserving functors (\cite[Def. 5.5.3.1]{HTT}), and the transition maps are the inverse image functors. For example, the motivic $\infty$-category $\mathcal{SH}$ is continuous (\cite[Prop. C.12(4)]{Hoy}). 

\subsubsection{}
Recall that a motivic $\infty$-catgegory $\mathcal{T}$ is \textbf{semi-separated} if for any finite, radicial and surjective morphism of schemes $f:T\to S$, the functor $f^*:\mathcal{T}(S)\to\mathcal{T}(T)$ is conservative (\cite[D\'ef. 2.1.160]{Ayo}, \cite[D\'ef. 1.1]{Ayo2} and \cite[Def. 2.1.7]{CD}). By \cite[Prop. 2.1.163]{Ayo}, this implies that for all morphisms $f$, the functor $f^*$ is an equivalence. 

\subsubsection{}
Recall that a morphism of schemes is a \emph{universal homeomorphism} if and only if it is integral, radicial and surjective (\cite[Cor. 18.12.11]{EGA4}). 
If $\mathcal{T}$ is continuous and semi-separated, then the functor $f^*$ is an equivalence for all universal homeomorphisms $f$, see \cite{EK}. 

\subsubsection{}
Recall the following topological invariance of $\mathcal{SH}$ due to Elmanto-Khan:
\begin{lemma}[\textrm{\cite[Th. 2.1.1]{EK}}]
\label{lm:EK211}
For $\mathcal{P}$ a collection of prime numbers, $\mathcal{SH}[\mathcal{P}^{-1}]$ is semi-separated when restricted to schemes where every prime outside $\mathcal{P}$ is invertible. 

In particular, for any prime $p$, $\mathcal{SH}[1/p]$ is semi-separated when restricted to schemes of characteristic $p$.
\end{lemma}

\section{The motivic nearby cycle}

\subsection{Specialization systems}
\label{sec:sps}
\subsubsection{}
Let $\mathcal{T}$ be a motivic $\infty$-catgegory over an $\infty$-category of schemes $\mathcal{S}$. In particular, $\mathcal{T}$ gives rise to a functor 
\begin{align}
\mathcal{T}^*:\mathcal{S}^{op}\to CAlg(Pr_{Stab})
\end{align}
which has a right adjoint $\mathcal{T}_*:\mathcal{S}\to CAlg(Pr_{Stab})$. In more concrete terms, for a scheme $S$ there is an associated presentable stable $\infty$-category $\mathcal{T}(S)$, and for a morphism $f:T\to S$ in $\mathcal{S}$, the functor $f^*:=\mathcal{T}^*(f):\mathcal{T}(S)\to\mathcal{T}(T)$ has a right adjoint $f_*:=\mathcal{T}_*(f):\mathcal{T}(T)\to\mathcal{T}(S)$.

\subsubsection{}
Recall our assumptions on the base scheme $B$ and the notations in~\ref{num:notB}.
Let $\mathcal{T}_1$ be a motivic $\infty$-catgegory over $Sch/\eta$ and $\mathcal{T}_2$ be a motivic $\infty$-catgegory over $Sch/\sigma$. We then have a functor
\begin{align}
j_*\mathcal{T}_1^*:
(Sch/B)^{op}
\xrightarrow{-\times_B\eta}
(Sch/\eta)^{op}
\xrightarrow{\mathcal{T}_1^*}
CAlg(Pr_{Stab})
\end{align}
where the first functor is given by the fiber product over $B$ by $\eta$. Similarly we have a functor
\begin{align}
i_*\mathcal{T}_2^*:
(Sch/B)^{op}
\xrightarrow{-\times_B\sigma}
(Sch/\sigma)^{op}
\xrightarrow{\mathcal{T}_2^*}
CAlg(Pr_{Stab})
\end{align}
where the first functor is given by the fiber product over $B$ by $\sigma$.

\subsubsection{}
Recall that if $\mathcal{B}$ and $\mathcal{C}$ are two $(\infty,1)$ categories and $F,G:\mathcal{B}\to\mathcal{C}$ are two functors, a \emph{lax natural transformation} $\eta:F\Rightarrow G$ is a functor $\eta:\mathcal{B}\to Fun(\Delta^1,\mathcal{C})$ such that $s\circ\eta=F$ and $t\circ\eta=G$, where $s,t:Fun(\Delta^1,\mathcal{C})\to\mathcal{C}$ assign the source and the target of a morphism, respectively (\cite[Definition 1.3]{JFS}).

We call a \emph{specialization presystem} from $\mathcal{T}_1$ to $\mathcal{T}_2$ a lax natural transformation $\psi:j_*\mathcal{T}_1^*\Rightarrow i_*\mathcal{T}_2^*$.

\subsubsection{}
In more concrete terms, a specialization presystem $\psi$ is the data of
\begin{itemize}
\item
for every morphism $X\xrightarrow{f}B$, an exact functor $\psi_f:\mathcal{T}_1(X_\eta)\to\mathcal{T}_2(X_\sigma)$
\item 
for every composable morphism $Y\xrightarrow{g}X\xrightarrow{f}B$, a natural transformation 
\begin{align}
\label{eq:psialpha}
g_\sigma^*\psi_f\to \psi_{f\circ g}g_\eta^*
\end{align}
compatible with compositions of $f$-morphisms in a homotopy coherent way.
\end{itemize}
By adjunction, the map~\eqref{eq:psialpha} gives rise to a map
\begin{align}
\label{eq:psibeta}
\psi_fg_{\eta*}\to g_{\sigma*}\psi_{f\circ g}.
\end{align}

\begin{definition}
A \textbf{specialization system} is a specialization presystem such that
\begin{enumerate}
\item
For every smooth morphism $g$, the map~\eqref{eq:psialpha} is invertible.
\item
For every proper morphism $g$, the map~\eqref{eq:psibeta} is invertible.
\end{enumerate}
A morphism of specialization systems $\psi\to\psi'$ is a morphism of lax natural transformations.
\end{definition}

\subsubsection{}
\label{num:spthom}
By \cite[Prop. 3.1.7]{Ayo}, specialization systems are compatible with the formation of Thom spaces of vector bundles: if $X$ is a $B$-scheme and $v$ is a virtual vector bundle over $X$, then there is a canonical isomorphism $Th(v_{|X_\sigma})\otimes\psi_X(-)\simeq\psi_X(-\otimes Th(v_{|X_\eta}))$. In particular, specialization systems commute with the formation of Tate twists.

\subsubsection{}
\label{num:combc}
Let $\mathcal{P}$ be a class of morphism of schemes. We say that be a specialization system $\psi$ commutes with $\mathcal{P}$-base changes if for any morphism $g$ in $\mathcal{P}$, the map~\eqref{eq:psialpha} is invertible.

\begin{lemma}
\label{lm:sptopinv}
Let $\psi$ be a specialization presystem from $\mathcal{T}_1$ to $\mathcal{T}_2$.
If both $\mathcal{T}_1$ and $\mathcal{T}_2$ are continuous and semi-separated, then $\psi$ commutes with base changes by universal homeomorphisms. 
\end{lemma}
\proof
Applying the functor $g_{\sigma*}$ (which is an equivalence of categories), the map $\alpha_g$ becomes
\begin{align}
\psi_f\simeq g_{\sigma*}g_\sigma^*\psi_f\to g_{\sigma*}\psi_{f\circ g}g_\eta^*\simeq \psi_fg_{\eta*}g_\eta^*
\end{align}
which agrees with the unit map of the adjunction $(g_\eta^*, g_{\eta*})$ by naturality of the transformation $g\mapsto\alpha_g$, and therefore is an isomorphism. 
\endproof

\subsubsection{}
\label{num:spst}
In what follows, we only consider the case where $\mathcal{T}_1$ and $\mathcal{T}_2$ arise in the following way: for $\mathcal{T}$ a motivic $\infty$-category over $Sch/B$, define
\begin{align}
\mathcal{T}_1=\mathcal{T}_{|\eta}:(Sch/\eta)^{op}\to (Sch/B)^{op}\to CAlg(Pr_{Stab})
\end{align}
\begin{align}
\mathcal{T}_2=\mathcal{T}_{|\sigma}:(Sch/\sigma)^{op}\to (Sch/B)^{op}\to CAlg(Pr_{Stab})
\end{align}
as the restrictions of $\mathcal{T}$ to $Sch/\eta$ and $Sch/\sigma$ respectively. We then use the terminology of specialization systems over $\mathcal{T}$.

\subsubsection{}
\label{num:cansp}
In the setting of~\ref{num:spst}, the functor
\begin{align}
\begin{split}
(Sch/B)^{op}&\to Fun(\Delta^1,CAlg(Pr_{Stab}))\\
X&\mapsto [\mathcal{T}(X_\eta)\xrightarrow{i_{X}^*j_{X*}}\mathcal{T}(X_\sigma)]
\end{split}
\end{align}
defines a specialization presystem
\begin{align}
\chi=i^*j_*:j_*\mathcal{T}_{|\eta}\Rightarrow i_*\mathcal{T}_{|\sigma}.
\end{align}
By virtue of smooth and proper base change properties of $\mathcal{T}$, $\chi$ is indeed a specialization system, called the \textbf{canonical specialization system}.

\subsection{Ayoub's construction}
\label{sec:ayocons}
\subsubsection{}
In this section, $B=\mathbb{A}^1_k$, where $k$ is a field of characteristic $p$. We recall Ayoub's construction of the motivic nearby cycle functor in \cite[\S3.4 and 3.5]{Ayo}.

\subsubsection{}
\label{num:cobar}
First recall Rector's \emph{geometric cobar construction} (\cite[2.3]{Rec}): if $\mathcal{C}$ is a category with finite products, 
$X\xrightarrow{f}B\xleftarrow{g}Y$ is a diagram in $\mathcal{C}$, then there is a cosimplicial object $X\widetilde{\times}_BY$ in $\mathcal{C}$ such that for $n\geqslant0$,
\begin{enumerate}
\item $(X\widetilde{\times}_BY)^n=X\times B\times\cdots\times B\times Y$, with $n$ copies of $B$.
\item For an object $T$ of $\mathcal{C}$, and maps $x:T\to X$, $y:T\to Y$, $b_1,\cdots,b_n:T\to B$ in $\mathcal{C}$, the cofaces $d^i:(X\widetilde{\times}_BY)^n\to (X\widetilde{\times}_BY)^{n+1}$ and codegeneracies $s^i:(X\widetilde{\times}_BY)^n\to (X\widetilde{\times}_BY)^{n-1}$ act as
\begin{align}
d^i(x,b_1,\cdots,b_n,y)=
\begin{cases}
(x,f(x),b_1,\cdots,b_n,y) & i=0 \\
(x,b_1,\cdots,b_i,b_i,\cdots,b_n,y) & 1\leqslant i\leqslant n\\
(x,b_1,\cdots,b_n,g(y),y) & i=n+1
\end{cases}
\end{align}
\begin{align}
s^i(x,b_1,\cdots,b_n,y)=(x,b_1,\cdots,b_i,b_{i+2},\cdots,b_n,y), 0\leqslant i\leqslant n-1.
\end{align}
\end{enumerate}

\subsubsection{}
For a scheme $S$, we have a diagram of smooth $\mathbb{G}_{m,S}$-schemes
\begin{align}
\label{eq:diaggm}
\mathbb{G}_{m,S}\xrightarrow{\delta}\mathbb{G}_{m,S}\times_S\mathbb{G}_{m,S}\xleftarrow{(id,1)}\mathbb{G}_{m,S}
\end{align}
where $\delta$ is the diagonal morphism, and the scheme $\mathbb{G}_{m,S}\times_S\mathbb{G}_{m,S}$ is seen as a smooth $\mathbb{G}_{m,S}$-scheme via the first projection. Applying the construction in~\ref{num:cobar} to the diagram~\eqref{eq:diaggm} and $\mathcal{C}=Sm/\mathbb{G}_{m,S}$ the category of smooth $\mathbb{G}_{m,S}$-schemes, we obtain a cosimplicial object $\mathcal{B}^\bullet\in Sm/\mathbb{G}_{m,S}^{\Delta}$, which enhances to a map between nerves
\begin{align}
N(\mathcal{B}^\bullet):N(\Delta)\to N(Sm/\mathbb{G}_{m,S}).
\end{align}
Denote by 
\begin{align}
\label{eq:infsus}
\Sigma^\infty_{\mathbb{G}_{m,S}}: 
N(Sm/\mathbb{G}_{m,S})
\to 
\mathcal{T}(Sm/\mathbb{G}_{m,S})
\end{align}
the infinite suspension map (\cite[Cor. 2.39]{Rob}). 
In terms of the six functors, the functor $\Sigma^\infty$ in~\eqref{eq:infsus} sends a smooth morphism $f:T\to \mathbb{G}_{m,S}$ to $f_\#\mathbbold{1}_T=f_!f^!\mathbbold{1}_T$, the ``motive'' of $T$ over $\mathbb{G}_{m,S}$, see \cite[1.1.34]{CD}.
Composing $N(\mathcal{B}^\bullet)$ with $\Sigma^\infty$ then gives a (homotopy coherent) cosimplicial object of $\mathcal{T}(\mathbb{G}_{m,S})$:
\begin{align}
\label{eq:sigmaA}
\mathcal{A}^\bullet=\Sigma^\infty N(\mathcal{B}^\bullet):
N(\Delta)\to\mathcal{T}(\mathbb{G}_{m,S}).
\end{align}

\subsubsection{}
\label{num:simpdual}
Recall that an object $K$ in a closed symmetric monoidal $\infty$-category $\mathcal{C}$ is \textbf{dualizable} if for any object $L\in\mathcal{C}$, the following canonical map is an isomorphism (see \cite[\S 4.6.1]{HA}):
\begin{align}
\underline{Hom}(K,\mathbbold{1}_X)\otimes L
\to 
\underline{Hom}(K,L).
\end{align}

For any scheme $T$, the infinite suspension of the smooth $T$-scheme $\mathbb{G}_{m,T}$ satisfies
\begin{align}
\Sigma^\infty_T\mathbb{G}_{m,T}\simeq \mathbbold{1}_T\oplus\mathbbold{1}_T(1)[1],
\end{align}
and therefore is a dualizable object in $\mathcal{T}(T)$. Consequently, the cosimplicial object $\mathcal{A}^\bullet$ in~\eqref{eq:sigmaA} factors through the full subcategory of dualizable objects in $\mathcal{T}(\mathbb{G}_{m,S})$.

\subsubsection{}
The \textbf{unipotent motivic nearby cycle functor} is the specialization system defined by
\begin{align}
\begin{split}
\Upsilon_X:\mathcal{T}(X_\eta)&\to\mathcal{T}(X_\sigma)\\
E&\mapsto\operatorname{colim}i_X^*j_{X*}\underline{Hom}(f_\eta^*\mathcal{A}^\bullet,E).
\end{split}
\end{align}

\subsubsection{}
The next construction takes into account finite \'etale covers, which reflect the discrete aspect of the motivic Galois group. Let $\mathbb{N}^\times$ be the category of positive integers, where there is a unique morphism $m\to n$ whenever $m | n$. Consider the functor $\mathcal{E}:\mathbb{N}^\times\to Sm/\mathbb{G}_{m,S}$, 
\begin{itemize}
\item
$\mathcal{E}$ sends every integer $n$ to the scheme $\mathbb{G}_{m,S}\times_S\mathbb{G}_{m,S}$, seen as a smooth $\mathbb{G}_{m,S}$-scheme via the first projection.
\item
$\mathcal{E}$ sends a morphism $m\to n$ to the morphism of schemes $\mathbb{G}_{m,S}\times_S\mathbb{G}_{m,S}\xrightarrow{(id,\frac{n}{m})}\mathbb{G}_{m,S}\times_S\mathbb{G}_{m,S}$ given by the identity on the first factor and $\frac{n}{m}$-th power map on the second factor.
\end{itemize}

\subsubsection{}
For a scheme $S$, we have a diagram in the category of functors $\mathbb{N}^\times\to Sm/\mathbb{G}_{m,S}$
\begin{align}
\label{eq:powergm}
\mathbb{G}_{m,S}\xrightarrow{\phi}\mathcal{E}\xleftarrow{(id,1)}\mathbb{G}_{m,S}.
\end{align}
Here the two terms $\mathbb{G}_{m,S}$ are considered as constant functors, and the map $\phi$ sends an integer $n$ to the morphism $\mathbb{G}_{m,S}\xrightarrow{(id,n)}\mathbb{G}_{m,S}\times_S\mathbb{G}_{m,S}$. Applying the construction in~\ref{num:cobar} to the diagram~\eqref{eq:powergm}, we obtain a functor $\mathcal{G}:\Delta\times\mathbb{N}^\times\to Sm/\mathbb{G}_{m,S}$. Similar to the construction in~\eqref{eq:sigmaA}, we define
\begin{align}
\mathcal{R}^\bullet=\Sigma^\infty N(\mathcal{G}^\bullet):
N(\Delta)\times N(\mathbb{N}^\times)=N(\Delta\times\mathbb{N}^\times)\to\mathcal{T}(\mathbb{G}_{m,S}).
\end{align}

\subsubsection{}
The \textbf{(tame) motivic nearby cycle functor} is the specialization system defined by
\begin{align}
\begin{split}
\Psi_X:\mathcal{T}(X_\eta)&\to\mathcal{T}(X_\sigma)\\
E&\mapsto\operatorname{colim}i_X^*j_{X*}\underline{Hom}(f_\eta^*\mathcal{R}^\bullet,E).
\end{split}
\end{align}

\subsubsection{}
The following description of the functor $\Psi$ is practical (\cite[Rem. 3.5.7]{Ayo}, \cite[Def. 1.4.11]{Ayo}).
For any integer $n\geqslant1$, denote by $e^n:\mathbb{A}^1_k\to \mathbb{A}^1_k$ the $n$-th power map. For a morphism $f:X\to \mathbb{A}^1_k$, form the Cartesian square
\begin{align}
\begin{split}
  \xymatrix@=10pt{
    X^n \ar[r]^-{e^n_X} \ar[d]_-{f^n} & X \ar[d]^-{f} \\
    \mathbb{A}^1_k \ar[r]^-{e^n} & \mathbb{A}^1_k.
  }
\end{split}
\end{align}
Then one can show that there is a canonical identification
\begin{align}
\Psi_X=\operatorname{colim}_{n}\Upsilon_{f^n}(e^n_X)_\eta^*.
\end{align}

\subsubsection{}
Ayoub shows that the three specialization systems $\chi$, $\Upsilon$ and $\Psi$ have a lax symmetric monoidal structure. For example, the lax monoidal structure of $\chi$ follows from the fact that for every scheme $X$, the functor $i_X^*$ is symmetric monoidal and the functor $j_{X*}$, being the right adjoint of a symmetric monoidal functor, is lax symmetric monoidal. For the case of $\Upsilon$ and $\Psi$, see (\cite[Prop. 3.2.17]{Ayo}. Furthermore, there are canonical isomorphisms (\cite[Prop. 3.4.9, Lemma 3.5.10]{Ayo})
\begin{align}
\label{eq:psiid1}
\mathbbold{1}_{\sigma}\simeq \Upsilon_{id}(\mathbbold{1}_{\eta}) \textrm{ and } \mathbbold{1}_{\sigma}\simeq \Psi_{id}(\mathbbold{1}_{\eta}).
\end{align}

\section{Local acyclicity and motivic nearby cycles}

\subsection{A K\"unneth formula}
\label{sec:kun}
\subsubsection{}
In characteristic $0$, Ayoub proves that the lax-monoidal functor $\Psi$ is indeed monoidal over the base scheme $\mathbb{A}^1$ (\cite[Cor. 3.5.18]{Ayo}). This property fails in positive characteristic, already for \'etale cohomology (see \cite[Rem. 2.2]{Ill}). The goal of this section is to prove a K\"unneth formula while restricting to certain motives with nice properties (see Theorem~\ref{th:lakun}).

\begin{definition}[\textrm{\cite[Def. 2.1.7]{JY}}]
\label{def:locacy}
Let $f:X\to S$ be a morphism of schemes of finite type and let $K\in \mathcal{T}(X)$. We say that $K$ is \textbf{strongly locally acyclic} (abbreviated as \textbf{SLA}) relatively to $f$ if for any morphism of finite type $p:T\to S$ with base change $q:X\times_ST\to X$, for any object $L\in\mathcal{T}(T)$, the canonical map
\begin{align}
\label{eq:lacan}
K\otimes (p_*L)_{|X}\to q_*(K\boxtimes_SL)
\end{align}
is an isomorphism. We say that $K$ is \textbf{universally strongly locally acyclic} (abbreviated as \textbf{USLA}) relatively to $f$ (or USLA over $S$) if for any morphism $T\to S$, the base change $K_{|X\times_ST}$ is strongly locally acyclic over $T$.
\end{definition}

\subsubsection{}
\label{num:lacond}
Here are two important cases where the USLA condition holds:
\begin{itemize}
\item
If $f:X\to S$ is smooth, any dualizable object in $\mathcal{T}(X)$ is USLA over $S$.
\item
If $S$ is the spectrum of a field $k$ of exponential characteristic $p$, and assume that either $k$ is a perfect field which satisfies strong resolution of singularities, or that $\mathcal{T}$ is $\mathbb{Z}[1/p]$-linear, then every object of $\mathcal{T}(X)$ is USLA over $k$.
\end{itemize}
Indeed, the first case follows from smooth base change and \cite[Prop. 3.2]{FHM}, and the second case is \cite[Cor. 2.1.14]{JY}.

\subsubsection{}
The next lemma shows that the USLA property is stable under exterior products:
\begin{lemma}
\label{lm:laprod}
Let $X$ and $Y$ be two $S$-schemes and let $K\in\mathcal{T}(X)$ and $L\in \mathcal{T}(Y)$ be two objects. If both $K$ and $L$ are USLA over $S$, then $K\boxtimes_SL\in\mathcal{T}(X\times_SY)$ is USLA over $S$.
\end{lemma}
\proof
Let $r:T\to S$ be a morphism of finite type, and let $M\in \mathcal{T}(T)$ be an object. Denote by $XY=X\times_SY$, and $s:X_T\to X$ and $t:(XY)_T\to XY$ the base changes of $r$. Then by the USLA property of $K$ and $L$ there are isomorphisms
\begin{align}
\label{eq:labox}
\begin{split}
(K\boxtimes_SL)\otimes(r_*M)_{|XY}
&\simeq 
L_{|XY}\otimes(K\otimes (r_*M)_{|X})_{|XY}
\simeq
L_{|XY}\otimes(s_*(K\boxtimes_SM))_{|XY}\\
&\simeq
t_*(L_{|XY}\boxtimes_X(K\boxtimes_SM))
\simeq
t_*((K\boxtimes_SL)\boxtimes_SM)
%_
\end{split}
\end{align}
and by naturality the map~\eqref{eq:labox} agrees with the canonical map~\eqref{eq:lacan}, which shows that $K\boxtimes_SL$ is strongly locally acyclic over $S$. The same argument shows that this property also holds after any base change over $S$, which shows that $K\boxtimes_SL$ is USLA over $S$.
\endproof

\subsubsection{}
Now we work over $B=\mathbb{A}^1_k$. For a $B$-scheme $X$ and $K\in\mathcal{T}(X)$, we denote $K_\eta=j_X^*K$ and $K_\sigma=i_X^*K$.
\begin{proposition}
\label{prop:lakun}
Let $X$ be an $B$-scheme and let $K\in\mathcal{T}(X)$ be an object. 
\begin{enumerate}
\item
If $K$ is SLA over $B$, then the following canonical map is an isomorphism:
\begin{align}
\label{eq:upsilonla}
K_{\sigma}\to\Upsilon_X(K_{\eta});
\end{align}
\item
\label{num:Psila}
If $K$ is USLA over $B$, then the following canonical map is an isomorphism:
\begin{align}
K_{\sigma}\to\Psi_X(K_{\eta}).
\end{align}
\end{enumerate}
\end{proposition}
\proof
\begin{enumerate}
\item
We have 
\begin{align}
\begin{split}
\Upsilon_X(K_{\eta})
&=
\operatorname{colim}i_X^*j_{X*}\underline{Hom}(f_\eta^*\mathcal{A}^\bullet,j_X^*K) \\
&\simeq
\operatorname{colim}i_X^*j_{X*}(j_X^*K\otimes f_\eta^*\underline{Hom}(\mathcal{A}^\bullet,\mathbbold{1}_{\eta}))\\
&\simeq
\operatorname{colim}i_X^*(K\otimes f^*j_*\underline{Hom}(\mathcal{A}^\bullet,\mathbbold{1}_{\eta}))\\
&\simeq
i_X^*K\otimes\operatorname{colim} f_\sigma^*i^*j_*\underline{Hom}(\mathcal{A}^\bullet,\mathbbold{1}_{\eta})\\
&\simeq
i_X^*K\otimes f_\sigma^*\operatorname{colim} i^*j_*\underline{Hom}(\mathcal{A}^\bullet,\mathbbold{1}_{\eta})\\
&=
i_X^*K\otimes f_\sigma^*\Upsilon_S\mathbbold{1}_\eta
\overset{\eqref{eq:psiid1}}{\simeq}
K_\sigma.
\end{split}
\end{align}
Here the first isomorphism follows from the fact that the cosimplicial object $\mathcal{A}^\bullet$ factors through dualizable objects (\ref{num:simpdual}), the second isomorphism follows from the SLA condition (Definition~\ref{def:locacy}), and the last two isomorphisms follow from the fact that the functors $i_X^*$,$\otimes$ and $f_\sigma^*$ have right adjoints and therefore commute with colimits.
\item
Using the USLA condition and the isomorphism~\eqref{eq:upsilonla}, we have 
\begin{align}
\begin{split}
\Psi_X(K_\eta)
&=
\operatorname{colim} \Upsilon_{f^n}(e^n_X)_\eta^*j_X^*K
=
\operatorname{colim} \Upsilon_{f^n}j_{X^n}^*(e^n)^*K\\
&\simeq
\operatorname{colim} i_{X^n}^*(e^n)^*K=i_X^*K=K_\sigma
\end{split}
\end{align}
which proves the result.
\end{enumerate}
\endproof

\begin{remark}
\label{rm:lavan}
Similar to the classical situation, one can define the \emph{motivic vanishing cycle functor} $\Phi_X:\mathcal{T}(X)\to\mathcal{T}(X_\sigma)$ as the cofiber of the map $i_X^*\to\Psi_Xj_X^*$. In \'etale cohomology, it is well-known that the stalks of the vanishing cycle complex are precisely the obstruction of being locally acyclic (\cite[Prop. 2.7]{Sai}, \cite[Cor. 6.6]{LZ}). In the motivic context, however, we do not know if the converse of Proposition~\ref{prop:lakun} is true, nor if it is reasonable to expect so.
\end{remark}

\begin{theorem}
\label{th:lakun}
Let $\psi$ be one of the specialization systems $\Upsilon$ or $\Psi$.
Let $X$ and $Y$ be two $B$-schemes and let $K\in \mathcal{T}(X)$ and $L\in \mathcal{T}(Y)$ be two objects. If both $K$ and $L$ are USLA over $B$, then the following canonical map induced by the lax-monoidal structure of $\psi$ is an isomorphism:
\begin{align}
\label{eq:psikunla}
\psi_X(K_\eta)\boxtimes_B\psi_Y(L_\eta)
\to
\psi_{X\times_BY}(K\boxtimes_BL)_\eta.
\end{align}
\end{theorem}

\proof
By Lemma~\ref{lm:laprod}, the object $K\boxtimes_BL$ is USLA over $B$. By Proposition~\ref{prop:lakun}, we have
\begin{align}
\label{eq:psiisola}
\psi_X(K_\eta)\boxtimes_B\psi_Y(L_\eta)
\simeq
K_\sigma\boxtimes_BL_\sigma
=
(K\boxtimes_BL)_\sigma
\simeq
\psi_{X\times_BY}(K\boxtimes_BL)_\eta.
\end{align}
and it is easy to see that the maps~\eqref{eq:psiisola} and~\eqref{eq:psikunla} agree, which proves the result.
\endproof

\subsubsection{}
In particular, for $X=Y=B$, we deduce from Theorem~\ref{th:lakun} the following
\begin{corollary}
\label{cor:psimon}
The functor $\Psi_{Id}j^*:\mathcal{T}(B)\to\mathcal{T}(\sigma)$ is symmetric monoidal when restricted to the full subcategory of objects that are USLA relatively to the identity morphism $id_B$.
\end{corollary}

\subsubsection{}
By~\ref{num:lacond}, any dualizable object in $\mathcal{T}(B)$ is USLA relatively to the identity morphism $id_B$. We will see later in Proposition~\ref{prop:ladual} that the converse is also true under resolution of singularities.

\subsection{More on local acyclicity}

\subsubsection{}
In this section we give more results on the local acyclicity properties.
We first study local acyclicity relative to the identity morphism. We begin with a lemma (see \cite[Prop. 6.3.5]{JY} for a similar statement):
\begin{lemma}
\label{lm:latr}
Let $X$ be a scheme and let $K\in\mathcal{T}(X)$ be USLA relatively to the identity morphism $id_X$. Then for any separated morphism of finite type $f:Y\to X$, the following canonical map is an isomorphism:
\begin{align}
f^*K\otimes f^!\mathbbold{1}_X
\to
f^!K
\end{align}
\end{lemma}
\proof
By working Zariski locally on $Y$, we may assume that the morphism $f$ is quasi-projective.
Choose a factorization of $f$ as $Y\xrightarrow{i}P\xrightarrow{g}X$, where $i$ is a closed immersion and $g$ is smooth. Let $j:U\to P$ be the open complement of $i$. Since $K$ is USLA over $X$, we have a canonical isomorphism $g^*K\otimes j_*\mathbbold{1}_U\simeq j_*j^*g^*K$. By localization we deduce a canonical isomorphism $g^*K\otimes i_*i^!\mathbbold{1}_P\simeq i_*i^!g^*K$, and applying the functor $i^*$ we obtain an isomorphism $i^*g^*K\otimes i^!\mathbbold{1}_P\simeq i^!g^*K$. On the other hand, since $g$ is smooth, the object $g^!\mathbbold{1}_X$ is dualizable, and by \cite[5.4]{FHM} we have
\begin{align}
\begin{split}
f^*K\otimes f^!\mathbbold{1}_X
&=
i^*g^*K\otimes i^!g^!\mathbbold{1}_X
\simeq
i^*g^*K\otimes i^!\mathbbold{1}_P\otimes i^*g^!\mathbbold{1}_X\\
&\simeq
i^!g^*K\otimes i^*g^!\mathbbold{1}_X
\simeq
i^!(g^*K\otimes g^!\mathbbold{1}_X)
\simeq
i^!g^!K
=
f^!K
\end{split}
\end{align}
which proves the result.
\endproof

\subsubsection{}
The following result is due to D.-C. Cisinski and the second-named author, and a more general version will appear in the forthcoming work~\cite{CY}. We include a proof for our case here for completeness:
\begin{proposition}
\label{prop:ladual}
Let $k$ be a field of exponential characteristic $p$ and $X$ be a regular $k$-scheme.
Assume that one of the following conditions hold:
\begin{enumerate}
\item
\label{num:dualres}
$\mathcal{T}$ is continuous and $k$ is a perfect field which satisfies strong resolution of singularities.
\item
\label{num:dualp}
$\mathcal{T}$ is $\mathbb{Z}[1/p]$-linear, and $X$ is smooth over a finite extension of $k$.
\end{enumerate}
Then an object $K\in\mathcal{T}(X)$ is USLA relatively to the identity morphism $id_X$ if and only if $K$ is a dualizable object.
\end{proposition}
\proof
If $K$ is dualizable then $K$ is USLA over $X$ by~\ref{num:lacond}.
For the converse, we first prove the case~\eqref{num:dualp}. Note that the unit object $\mathbbold{1}_X$ is a dualizing object in $\mathcal{T}(X)$. By \cite[Lemma 2.3.2]{Jin}, it suffices to show that for any projective morphism $f:Y\to X$ with $Y$ smooth over a finite extension of $k$, the following canonical map is an isomorphism:
\begin{align}
\underline{Hom}(K,\mathbbold{1}_X)\otimes f_*\mathbbold{1}_Y
\to 
\underline{Hom}(K,f_*\mathbbold{1}_Y).
\end{align}

On the one hand, by projection formula we have $\underline{Hom}(K,\mathbbold{1}_X)\otimes f_*\mathbbold{1}_Y\simeq f_*f^*\underline{Hom}(K,\mathbbold{1}_X)$.

On the other hand, by Lemma~\ref{lm:latr} we have $f^*K\otimes f^!\mathbbold{1}_X\simeq f^!K$. Since both $\mathbbold{1}_X$ and $f^!\mathbbold{1}_X$ are at the same time dualizing and $\otimes$-invertible objects, we have
\begin{align}
\begin{split}
\underline{Hom}(K,f_*\mathbbold{1}_Y)
&\simeq
f_*\underline{Hom}(f^*K,\mathbbold{1}_Y)
\simeq
f_*\underline{Hom}(f^*K\otimes f^!\mathbbold{1}_X,f^!\mathbbold{1}_X)\\
&\simeq
f_*\underline{Hom}(f^!K,f^!\mathbbold{1}_X)
\simeq
\underline{Hom}(f_!f^!K,\mathbbold{1}_X)
\simeq
f_*f^*\underline{Hom}(K,\mathbbold{1}_X)
\end{split}
\end{align}
which proves the result.

In the case~\eqref{num:dualres}, by Popescu's theorem (\cite[Th. 1.1]{Spi}) and continuity, we may assume that $X$ is smooth over $k$. Then the proof follows the same lines as the case~\eqref{num:dualp}, where we apply \cite[Prop. 2.2.27]{Ayo} instead of \cite[Lemma 2.3.2]{Jin}.%Lemma~\ref{lm:CD7.2}.
\endproof

\subsubsection{}
We now discuss further base change properties of motivic nearby cycles. 
Recall from~\ref{num:combc} that a specialization system commutes with base changes for a class of morphisms of schemes if the map~\eqref{eq:psialpha} is invertible.
By definition, any specialization system commutes with smooth base changes, and we have seen in Lemma~\ref{lm:sptopinv} that if the motivic $\infty$-category is semi-separated, then any specialization system commutes with base changes by finite, radicial and surjective morphisms. 
\begin{definition}
We say that a morphism of schemes $f:X\to S$ is \textbf{SLA} (resp. \textbf{USLA}) if the unit object is SLA (resp. USLA) relatively to $f$.
\end{definition}

\subsubsection{}
\label{num:labc}
To say that a morphism $f:X\to S$ is SLA amounts to say that for any Cartesian diagram of schemes
\begin{align}
\begin{split}
  \xymatrix@=10pt{
    Y \ar[r]^-{q} \ar[d]_-{g} & X \ar[d]^-{f} \\
    T \ar[r]^-{p} & S
  }
\end{split}
\end{align}
the canonical map $f^*p_*\to q_*g^*$ is invertible (which can be expressed as the corresponding diagram of pull-back functors being \emph{right adjointable} as in \cite[Def. 4.7.4.13]{HA}, or equivalently that $\mathcal{T}$ satisfies \emph{$f$-base change} as in \cite[Def. 2.1.3]{JY}). The USLA property further requires this property to hold after any base change over $S$.

\subsubsection{}
By~\ref{num:lacond}, any smooth morphism, or any morphism of the form $X\to k$, where either $k$ is a perfect field which satisfies strong resolution of singularities or $\mathcal{T}$ is $\mathbb{Z}[1/p]$-linear, is USLA.

\begin{proposition}
\label{prop:labc}
\begin{enumerate}
\item
The specialization systems $\chi$ and $\Upsilon$ commute with SLA base changes.
\item
The specialization system $\Psi$ commutes with USLA base changes.
\end{enumerate}
\end{proposition}
\proof
Let $Y\xrightarrow{g}X\xrightarrow{f}B$ be two composable morphisms. Assume that the morphism $g$ is SLA. By~\ref{num:labc}, we have
\begin{align}
\label{eq:SLAchi}
g_\sigma^*\chi_f=g_\sigma^*i_X^*j_{X*}=i_Y^*g^*j_{X*}\simeq i_Y^*j_{Y*}g_\eta^*=\chi_{f\circ g}g_\eta^*.
\end{align}
so the specialization system $\chi$ commutes with $g$-base change. Similarly we have
\begin{align}
\begin{split}
g_\sigma^*\Upsilon_f
&=
g_\sigma^*\operatorname{colim}\chi_f\underline{Hom}(f_\eta^*\mathcal{A}^\bullet,-)\\
&\simeq
\operatorname{colim}\chi_{(f\circ g)}g_\eta^*\underline{Hom}(f_\eta^*\mathcal{A}^\bullet,-)\\
&\simeq
\operatorname{colim}\chi_{f\circ g}\underline{Hom}((f\circ g)_\eta^*\mathcal{A}^\bullet,g_\eta^*(-))
=
\Upsilon_{f\circ g}g_\eta^*
\end{split}
\end{align}
where the first isomorphism follows from the fact that the functor $g^*$ commutes with $\chi$~\eqref{eq:SLAchi} and colimits, and the second isomorphism follows from the fact that the cosimplicial object $\mathcal{A}^\bullet$ factors through dualizable objects (\ref{num:simpdual}). The case of the specialization system $\Psi$ is similar and is left as exercise (see the proof of Proposition~\ref{prop:lakun}~\eqref{num:Psila}).
\endproof


\begin{thebibliography}{}

\bibitem[Ayo07a]{Ayo}
J. Ayoub,
\emph{Les six op\'erations de Grothendieck et le formalisme des cycles \'evanescents dans le monde motivique}, Ast\'erisque No. \textbf{314}, \textbf{315} (2007).

\bibitem[Ayo07b]{Ayo1}
J. Ayoub,
\emph{The motivic vanishing cycles and the conservation conjecture}, 
in \emph{Algebraic cycles and motives}, Volume 1. Selected papers of the EAGER conference, Leiden, Netherlands, August 30–September 3, 2004 on the occasion of the 75th birthday of Professor J. P. Murre. London Mathematical Society Lecture Note Series \textbf{343}, 3-54 (2007). 

\bibitem[Ayo14]{Ayo2}
J. Ayoub,
\emph{La r\'ealisation \'etale et les op\'erations de Grothendieck}, 
Ann. Sci. \'Ec. Norm. Sup\'er. (4) \textbf{47} (2014), no. 1, 1-145.

\bibitem[Ayo22]{Ayo3}
J. Ayoub,
\emph{Anabelian presentation of the motivic Galois group in characteristic zero}, 
preprint, 2022.

\bibitem[AIS17]{AIS}
J. Ayoub, F. Ivorra, J. Sebag, 
\emph{Motives of rigid analytic tubes and nearby motivic sheaves}, 
Ann. Sci. \'Ec. Norm. Sup\'er. (4) \textbf{50} (2017), no. 6, 1335-1382. 

\bibitem[Azo21]{Azo}
R. Azouri, 
\emph{Motivic Euler characteristic of nearby cycles and a generalized quadratic conductor formula}, 
arxiv preprint \url{https://arxiv.org/abs/2101.02686}.

\bibitem[CD19]{CD}
D.-C. Cisinski, F. D\'eglise,
\emph{Triangulated categories of motives}, 
Springer Monographs in Mathematics. Springer, Cham, 2019.

\bibitem[CY22]{CY}
D.-C. Cisinski, E. Yang,
\emph{Motivic singular support}, 
in preparation, 2022.

\bibitem[DFJK21]{DFJK}
F. D\'eglise, J. Fasel, F. Jin, A. Khan, 
\emph{On the rational motivic homotopy category}, 
J. Ec. Polytech. Math. \textbf{8} (2021), 533-583.

\bibitem[DG20]{DG}
B. Drew, M. Gallauer,
\emph{The Universal Six-Functor Formalism}, 
arXiv:2009.13610.

\bibitem[EGA4]{EGA4}
A. Grothendieck, 
\emph{\'El\'ements de g\'eom\'etrie alg\'ebrique. IV. \'Etude locale des sch\'emas et des morphismes de sch\'emas}, 
r\'edig\'es avec la collaboration de Jean Dieudonn\'e,
Inst. Hautes \'Etudes Sci. Publ. Math. No. \textbf{20}, \textbf{24}, \textbf{28}, \textbf{32}, 1964-1967.

\bibitem[EK20]{EK}
E. Elmanto, A. Khan,
\emph{Perfection in motivic homotopy theory},  Proc. Lond. Math. Soc. \textbf{120} (2020), no. 1, 28-38.

\bibitem[FHM03]{FHM}
H. Fausk, P. Hu, J. P. May,
\emph{Isomorphisms between left and right adjoints}, Theory Appl. Categ. \textbf{11} (2003), No. 4, 107-131.

\bibitem[Hoy15]{Hoy}
M. Hoyois,
\emph{A quadratic refinement of the Grothendieck-Lefschetz-Verdier trace formula}, Algebraic \& Geometric Topology \textbf{14} (2015), no.~6, 3603-3658.

\bibitem[Ill17]{Ill}
L. Illusie,
\emph{Around the Thom-Sebastiani theorem}, with an appendix by Weizhe Zheng,
Manuscripta Math. \textbf{152} (2017), no. 1-2, 61-125. 

\bibitem[Jin20]{Jin}
F. Jin,
\emph{Trace maps in motivic homotopy and local terms}, 
\href{https://arxiv.org/abs/2010.09292}{arXiv:2010.09292}.

\bibitem[JY21]{JY}
F. Jin, E. Yang,
\emph{K\"unneth formulas for motives and additivity of traces}, 
Adv. Math. \textbf{376} (2021), Article ID 107446.

\bibitem[JFS17]{JFS}
T. Johnson-Freyd, C. Scheimbauer,
\emph{(Op)lax natural transformations, twisted quantum field theories, and "even higher'' Morita categories}, 
Adv. Math. \textbf{307} (2017), 147-223. 

\bibitem[Kha16]{Kha}
A. Khan,
\emph{Motivic homotopy theory in derived algebraic geometry}, 
Ph.D. thesis, Universit\"at Duisburg-Essen, 2016, available at \url{https://www.preschema.com/thesis/}.

\bibitem[LPS21]{LPS}
M. Levine, S. Pepin-Lehalleur, V. Srinivas,
\emph{Euler characteristics of homogeneous and weighted-homogeneous hypersurfaces}, 
arxiv preprint \url{https://arxiv.org/abs/2101.00482}.

\bibitem[LZ19]{LZ}
Q. Lu, W. Zheng,
\emph{Duality and nearby cycles over general bases}, 
Duke Math. J. \textbf{168} (2019), no. 16, 3135-3213.

\bibitem[Lur09]{HTT}
J. Lurie,
\emph{Higher topos theory}, 
Annals of Mathematics Studies, \textbf{170}. Princeton University Press, Princeton, NJ, 2009.

\bibitem[Lur17]{HA}
J. Lurie,
\emph{Higher Algebra}, 
available at \url{https://www.math.ias.edu/~lurie/papers/HA.pdf}.

\bibitem[RZ82]{RZ}
M. Rapoport, T. Zink,
\emph{\"Uber die lokale Zetafunktion von Shimuravariet\"aten. Monodromiefiltration und verschwindende Zyklen in ungleicher Charakteristik}, 
Invent. Math. \textbf{68} (1982), no. 1, 21-101. 

\bibitem[Rec70]{Rec}
D. Rector, 
\emph{Steenrod operations in the Eilenberg-Moore spectral sequence}, 
Comment. Math. Helv. \textbf{45} (1970), 540-552.

\bibitem[Rob15]{Rob}
M. Robalo, 
\emph{$K$-theory and the bridge from motives to noncommutative motives}, 
Adv. Math. \textbf{269} (2015), 399550. 

\bibitem[Sai17]{Sai}
T. Saito,
\emph{The characteristic cycle and the singular support of a constructible sheaf}, Invent. Math. \textbf{207} (2017), no. 2, 597-695.

\bibitem[SGA7]{SGA7}
\emph{Groupes de monodromie en g\'eom\'etrie alg\'ebrique I, II}, S\'eminaire de G\'eom\'etrie Alg\'ebrique du Bois-Marie 1967-1969 (SGA 7). Dirigé par A. Grothendieck, P. Deligne et N. Katz. Avec la collaboration de M. Raynaud et D. S. Rim. Lecture Notes in Mathematics, Vol. \textbf{288}, \textbf{340}. Springer-Verlag, Berlin-New York, 1972-1973.

\bibitem[Spi99]{Spi}
M. Spivakovsky,
\emph{A new proof of D. Popescu’s theorem on smoothing of ring homomorphisms}, J. Amer. Math. Soc. \textbf{12} (1999), no. 2, p. 381-444.

\end{thebibliography}
\end{document}